\documentclass[12pt]{article}
\usepackage{epic,eepic,epsf,epsfig}
\usepackage{amsmath}
\usepackage{amssymb}
\usepackage{amsbsy}

\usepackage{amsfonts}
\newtheorem{theorem}{Theorem}%
\newtheorem{cor}[theorem]{Corollary}%
\newtheorem{remark}[theorem]{Remark}%

 \def\Omega{\Omega}

\def\f{\noindent}

\def\Aut{\hbox{\rm Aut}}

\def\Syl{\hbox{\rm Syl}}

\newcommand{\qed}{\mbox{\raisebox{0.7ex}{\fbox{}}} \vspace{4truemm}}
\def\demo{\f {\bf Proof.}\hskip10pt}

\begin{document}

\baselineskip 16pt

\title{\vspace{-1.2cm}
Finite groups in which every self-centralizing subgroup is a TI-subgroup or subnormal or has $p'$-order
\thanks{\scriptsize
This work was supported by Shandong Provincial Natural Science Foundation, China (ZR2017MA022 and ZR2020MA044)
and NSFC (11761079).
\newline
\hspace*{0.5cm} \scriptsize{E-mail address: jiangtaoshi@126.com\,(J. Shi).}}}

\author{Jiangtao Shi\\
\\
{\small{\em School of Mathematics and Information Sciences, Yantai University, Yantai 264005, P.R. China}}}
\date{ }

\maketitle \vspace{-.8cm}

\begin{abstract}
We first give complete characterizations of the structure of finite group $G$ in which every subgroup (or non-nilpotent subgroup,
or non-abelian subgroup) is a TI-subgroup or subnormal or has $p'$-order for a fixed prime divisor $p$ of $|G|$. Furthermore,
we prove that every self-centralizing subgroup (or non-nilpotent subgroup, or non-abelian subgroup) of $G$ is a TI-subgroup
or subnormal or has $p'$-order for a fixed prime divisor $p$ of $|G|$ if and only if every subgroup
(or non-nilpotent subgroup, or non-abelian subgroup) of $G$ is a TI-subgroup or subnormal or has $p'$-order. Based on these results,
we obtain the structure of finite group $G$ in which every self-centralizing
subgroup (or non-nilpotent subgroup, or non-abelian subgroup) is a TI-subgroup or subnormal or has $p'$-order for a fixed prime divisor
$p$ of $|G|$.

\medskip \f {\bf Keywords:} self-centralizing; TI-subgroup; subnormal; $p'$-order; non-nilpotent subgroup; non-abelian subgroup
\end{abstract}

\section{Introduction}

Throughout this paper all groups are assumed to be finite. Suppose that $G$ is a group and $H$ a subgroup of $G$, then $H$ is
termed to be a TI-subgroup of $G$ if $H^g\cap H=1$ or $H$ for each $g\in G$. It is clear that the TI-subgroup and the subnormal subgroup are two
relatively independent concepts in group theory. For TI-subgroups, Walls {\rm\cite{walls} classified groups in which
every subgroup is a TI-subgroup. As a generalization, Shi and Zhang {\rm\cite[Theorem 2]{s2}} characterized groups of even
order in which every subgroup of even order is a TI-subgroup. For subnormal subgroups, {\rm\cite[Theorem 5.2.4]{rob}} indicated
that a group $G$ is nilpotent if and only if every subgroup of $G$ is subnormal. In {\rm\cite{eb}} Ebert and Bauman characterized
groups in which every subgroup is subnormal or abnormal. Kurdachenko and Smith {\rm\cite{ku}} investigated general groups in which
every subgroup is either subnormal or self-normalizing. In {\rm\cite[Theorem 1]{ba}} Ballester-Bolinches
and Cossey described groups in which every subgroup is supersolvable or subnormal.

Combined the TI-property and the subnormality of subgroups together Shi and Zhang {\rm\cite[Theorem 1]{s1}} gave a complete
characterization of groups in which every subgroup is a TI-subgroup or subnormal. Combined the nilpotence, the normality and
the order of subgroups together Shi, Li and Shen {\rm\cite[Theorems 1.3, 1.4 and Theorem 1.7]{s5}}
investigated group $G$ in which every maximal subgroup is nilpotent or normal or has $p'$-order for a fixed prime divisor $p$
of $|G|$. Furthermore, Shi {\rm\cite[Theorem 1.1]{s6}} obtained a complete characterization of group $G$ in which every maximal
subgroup is nilpotent or a TI-subgroup or has $p'$-order for a fixed prime divisor $p$ of $|G|$.

In this paper, motivated by above researches, combining the TI-property, the subnormality and the order of subgroups together
we have the following result whose proof is given in Section~\ref{th01}.

\begin{theorem}\ \ \label{th1} Suppose that $G$ is a group and $p$ a fixed prime divisor of $|G|$. Then every subgroup of $G$
is a TI-subgroup or subnormal or has $p'$-order if and only if one of the following statements holds:

$(1)$ every subgroup of $G$ of order divisible by $p$ is subnormal in $G$;

$(2)$ $p=2$, $G=Z_q\rtimes\langle a\rangle$ is a Frobenius group with kernel $Z_q$ and complement $\langle a\rangle$,
where $q$ is an odd prime and $o(a)$ is an even number;

$(3)$ $p>2$, $G={Z_q}^r\rtimes(P\times H)$ is a Frobenius group with kernel ${Z_q}^r$ and complement $P\times H$,
where $q\neq p$ and $r\geq 1$, $P\in\Syl_p(G)$ and $P$ is cyclic, $H$ is either a cyclic group or a direct product
of a quaternion group $Q_8$ and a cyclic group of odd order, and every non-identity subgroup of $P$ acts irreducibly
on ${Z_q}^r$;

$(4)$ $p>2$, $G={Z_q}^r\rtimes(Z_p\rtimes H)$ is a Frobenius group with kernel ${Z_q}^r$ and complement $Z_p\rtimes H$,
where $q\neq p$ and $r>1$, $Z_p\in\Syl_p(G)$, $H$ is either a cyclic group or a direct product of a quaternion group
$Q_8$ and a cyclic group of odd order such that $[Z_p,H]\neq 1$, and $Z_p$ acts irreducibly on ${Z_q}^r$.
\end{theorem}

In {\rm\cite[Theorem 1 and Corollary 2]{s3}} Shi characterized groups in which every non-abelian subgroup is a TI-subgroup
or subnormal. In this paper, as a further generalization and extension, assume that every non-nilpotent subgroup (or non-abelian subgroup) of a group $G$ is a
TI-subgroup or subnormal or has $p'$-order for a fixed prime divisor $p$ of $|G|$, arguing as in the proof of
Theorem~\ref{th1}, we can obtain the following Theorem~\ref{th2} and Theorem~\ref{th3}, here we omit their proofs.

\begin{theorem}\ \ \label{th2} Suppose that $G$ is a group and $p$ a fixed prime divisor of $|G|$. Then every
non-nilpotent subgroup of $G$ is a TI-subgroup or subnormal or has $p'$-order if and only if one
of the following statements holds:

$(1)$ every non-nilpotent subgroup of $G$ of order divisible by $p$ is subnormal in $G$;

$(2)$ $p>2$, $G={Z_q}^r\rtimes(Z_p\rtimes H)$ is a Frobenius group with kernel ${Z_q}^r$ and complement $Z_p\rtimes H$,
where $q\neq p$ and $r>1$, $Z_p\in\Syl_p(G)$, $H$ is either a cyclic group or a direct product of a quaternion group
$Q_8$ and a cyclic group of odd order such that $[Z_p,H]\neq 1$, and $Z_p$ acts irreducibly on ${Z_q}^r$.
\end{theorem}

The following corollary is a direct consequence of Theorem~\ref{th2}.

\begin{cor}\ \ \label{c1} Suppose that $G$ is a group and $p$ the smallest prime divisor of $|G|$. Then every
non-nilpotent subgroup of $G$ is a TI-subgroup or subnormal or has $p'$-order if and only if every non-nilpotent
subgroup of $G$ of order divisible by $p$ is subnormal in $G$.
\end{cor}

\begin{theorem}\ \ \label{th3} Suppose that $G$ is a group and $p$ a fixed prime divisor of $|G|$. Then every
non-abelian subgroup of $G$ is a TI-subgroup or subnormal or has $p'$-order if and only if one
of the following statements holds:

$(1)$ every non-abelian subgroup of $G$ of order divisible by $p$ is subnormal in $G$;

$(2)$ $p>2$, $G={Z_q}^r\rtimes(P\times H)$ is a Frobenius group with kernel ${Z_q}^r$ and complement
$P\times H$, where $q\neq p$ and $r>1$, $P\in\Syl_p(G)$ and $P$ is cyclic, $H$ is a direct product of a quaternion
group $Q_8$ and a cyclic group of odd order, and every non-identity subgroup of $P$ acts irreducibly on ${Z_q}^r$;

$(3)$ $p>2$, $G={Z_q}^r\rtimes(Z_p\rtimes H)$ is a Frobenius group with kernel ${Z_q}^r$ and complement $Z_p\rtimes H$,
where $q\neq p$ and $r>1$, $Z_p\in\Syl_p(G)$, $H$ is either a cyclic group or a direct product of a quaternion group
$Q_8$ and a cyclic group of odd order such that $[Z_p,H]\neq 1$, and $Z_p$ acts irreducibly on ${Z_q}^r$.
\end{theorem}

Let $G$ be a group and $H$ a subgroup of $G$, then $H$ is said to be self-centralizing in $G$ if $C_G(H)\leq H$.
As an extension of {\rm\cite[Theorem 1 and Corollary 2]{s3}}, Sun, Lu and Meng {\rm\cite[Theorem 1.1]{slm}} proved that
if every self-centralizing non-abelian subgroup of a group $G$ is a TI-subgroup or subnormal then every
non-abelian subgroup of $G$ is subnormal. Furthermore, Shi and Li {\rm\cite[Theorem 1 and Theorem 2]{s4}} investigated
groups in which every self-centralizing non-nilpotent subgroup is a TI-subgroup or subnormal.

According to above results,
it is natural and interesting to characterize the structure of group $G$ in which every self-centralizing subgroup
(or non-nilpotent subgroup, or non-abelian subgroup) is a TI-subgroup or
subnormal or has $p'$-order for a fixed prime divisor $p$ of $|G|$. In this paper, we obtain the following Theorem~\ref{th001} which
can indicate the equivalent relation between group $G$ in which every self-centralizing subgroup is a TI-subgroup or
subnormal or has $p'$-order for a fixed prime divisor $p$ of $|G|$ and group $G$ in which every subgroup is a TI-subgroup or
subnormal or has $p'$-order for a fixed prime divisor $p$ of $|G|$.

\begin{theorem}\ \ \label{th001} Suppose that $G$ is a group and $p$ a fixed prime divisor of $|G|$. Then every self-centralizing subgroup of
$G$ is a TI-subgroup or subnormal or has $p'$-order if and only if every subgroup of $G$ is a TI-subgroup or subnormal or has $p'$-order.
\end{theorem}

The proof of Theorem~\ref{th001} is given in Section~\ref{th002}.

Arguing as in proof of Theorem~\ref{th001}, we can also obtain the following two results, here we omit their proofs.

\begin{theorem}\ \ \label{th012} Suppose that $G$ is a group and $p$ a fixed prime divisor of $|G|$. Then every self-centralizing non-nilpotent
subgroup of $G$ is a TI-subgroup or subnormal or has $p'$-order if and only if every non-nilpotent subgroup of $G$ is a TI-subgroup or subnormal
or has $p'$-order.
\end{theorem}

\begin{theorem}\ \ \label{th013} Suppose that $G$ is a group and $p$ a fixed prime divisor of $|G|$. Then every self-centralizing non-abelian
subgroup of $G$ is a TI-subgroup or subnormal or has $p'$-order if and only if every non-abelian subgroup of $G$ is a TI-subgroup or subnormal or has
$p'$-order.
\end{theorem}

\begin{remark}\ \ {\rm Theorems~\ref{th1}, \ref{th2}, \ref{th3}, \ref{th001}, \ref{th012} and \ref{th013} give the structure of group $G$ in which every self-centralizing
subgroup (or non-nilpotent subgroup, or non-abelian subgroup) is a TI-subgroup or subnormal or has $p'$-order for a fixed prime divisor
$p$ of $|G|$.}
\end{remark}

\section{Proof of Theorem~\ref{th1}}\label{th01}

\demo We first prove the necessity part.

Assume that $G$ has at least one subgroup of order divisible by $p$ that is not subnormal in $G$. Let $M$ be the largest
subgroup of $G$ of order divisible by $p$ that is not subnormal in $G$, one must have $M=N_G(M)$. Since $M$ is a TI-subgroup
of $G$ by the hypothesis, we get that $G$ is a Frobenius group with complement $M$. Let $N$ be the kernel, then $G=N\rtimes M$.

Let $P\in\Syl_p(M)$, obviously $P\in\Syl_p(G)$. Note that the order of the subgroup $N\rtimes P$ is divisible by $p$. By the
hypothesis, $N\rtimes P$ is a TI-subgroup of $G$ or subnormal in $G$. For the case when $N\rtimes P$ is a TI-subgroup of
$G$, one has $N\rtimes P\trianglelefteq G$ since $(N\rtimes P)^g\cap(N\rtimes P)=(N^g\rtimes P^g)\cap(N\rtimes P)=(N\rtimes P^g)\cap(N\rtimes P)
\geq N\neq 1$ for each $g\in G$. It follows that $P=(N\cap M)P=(N\rtimes P)\cap M\trianglelefteq M$. For another case when $N\rtimes P$
is subnormal in $G$, one has that $P$ is subnormal in $M$ and then $P\trianglelefteq M$ since $P\in\Syl_p(M)$. Thus we always
have $P\trianglelefteq M$. By Schur-Zassenhaus theorem (see {\rm\cite[Theorem 9.1.2]{rob}}, $M$ has a subgroup $H$ such that $M=PH$
and $P\cap H=1$, that is $M=P\rtimes H$.

{\bf Claim 1}: $M$ is maximal in $G$. It is clear that $M\leq N_G(P)<G$. Let $K$ be a maximal subgroup of $G$ such that
$M\leq N_G(P)\leq K$. If $M<K$, by the choice of $M$, one has that $K$ is subnormal in $G$. Then $K\trianglelefteq G$.
By Frattini-argument, one has $P\trianglelefteq G$, a contradiction. Thus $M=N_G(P)=K$ is maximal in $G$.

By Claim 1, $N$ is a minimal normal subgroup of $G$. Since $N$ is nilpotent by {\rm\cite[Theorem 10.5.6($i$)]{rob}},
one has that $N$ is an elementary abelian group. Assume $N={Z_q}^r$, where
$q\neq p$ and $r\geq 1$.

{\bf Claim 2}: $H$ is nilpotent. If $H=1$, $H$ is obvious nilpotent. Next assume $H>1$. For any maximal subgroup $H_1$ of $H$,
$N\rtimes(P\rtimes H_1)$ is a maximal subgroup of $G$ of order divisible by $p$. By the hypothesis, $N\rtimes(P\rtimes
H_1)$ is a TI-subgroup of $G$ or subnormal in $G$. Arguing as above, we can get $H_1\trianglelefteq H$.
By the choice of $H_1$, one has that $H$ is nilpotent.

{\bf Claim 3}: By conjugation every non-identity subgroup of $P$ acts irreducibly on $N$. Otherwise, assume that $P_1$
is a non-identity subgroup of $P$ and $N_1$ is a non-trivial subgroup of $N$ such that $P_1$ normalizes $N_1$.
Then $N_1\rtimes P_1$ is a subgroup of $G$ of order divisible by $p$. By the hypothesis, $N_1\rtimes P_1$ is a TI-subgroup of $G$ or
subnormal in $G$. One has that $N_1\rtimes P_1$ is also a TI-subgroup of $N\rtimes P_1$ or subnormal in $N\rtimes P_1$.
For the case when $N_1\rtimes P_1$ is a TI-subgroup of $N\rtimes P_1$. Note that $N_1\trianglelefteq N\rtimes P_1$ since $N$ is abelian.
Then $(N_1\rtimes P_1)^g\cap(N_1\rtimes P_1)=(N_1^g\rtimes P_1^g)\cap (N_1\rtimes P_1)=(N_1\rtimes P_1^g)\cap(N_1\rtimes P_1)
\geq N_1\neq 1$ for each $g\in N\rtimes P_1$. It follows that $N_1\rtimes P_1\trianglelefteq N\rtimes P_1$, one has either
$N_1\rtimes P_1\leq N$ or $N<N_1\rtimes P_1$. It is obvious that both of them are impossible. For another case when
$N_1\rtimes P_1$ is subnormal in $N\rtimes P_1$. Assume that $N_1\rtimes P_1=L_0\vartriangleleft L_1\vartriangleleft\cdots
\vartriangleleft L_s=N\rtimes P_1$ is a subnormal subgroups series, where $s\geq 1$. Observing that $N\nleq L_0$ but
$N\leq L_s$. Let $m$ be the smallest number from $1$ to $s$ such that $N\leq L_m$ but $N\nleq L_{m-1}$. Since $L_0\leq L_m$
but $L_0\nleq N$, one has $N<L_m$. Then $L_m=N\rtimes (L_m\cap P_1)$ is also a Frobenius group with kernel $N$ and complement
$L_m\cap P_1>1$. Since $L_{m-1}\trianglelefteq L_m$, one has either $L_{m-1}\leq N$ or $N<L_{m-1}$. It is also obvious that
both of them are impossible. Hence every non-identity subgroup of $P$ acts irreducibly on $N$.

{\bf Case I}: Assume $p=2$. Take an element $e$ of order $2$ in $P$. Since $C_G(N)=N$ as $G$ being a Frobenius group,
$e$ induces a fixed-point-free automorphism of $N$ of order 2. One has $g^e=g^{-1}$ for each non-identity element $g\in N$.
Note that every non-identity subgroup of $P$ acts irreducibly on $N$. It follows that $N=Z_q$. By N/C-theorem,
$M\cong G/N=N_G(N)/C_G(N)\lesssim\Aut(N)$ is a cyclic group. One has $G=Z_q\rtimes\langle a\rangle$, where $o(a)$
is an even number.

{\bf Case II}: Assume $p>2$.

{\bf Subcase (1)}: Assume that the order of every maximal subgroup of $M$ is divisible by $p$. Arguing as above, one
can get that every maximal subgroup of $M$ is normal in $M$. Then $M$ is nilpotent. One has $M=P\times H$ and then
$G={Z_q}^r\rtimes(P\times H)$. Since $p>2$, $P$ is cyclic by {\rm\cite[Theorem 10.5.6($ii)$]{rob}}. Let $M_0$ be
any non-trivial subgroup of $P\times H$ of order divisible by $p$, then $M_0=(M_0\cap P)\times(M_0\cap H)$, where
$M_0\cap P>1$. By the hypothesis, $M_0$ is a TI-subgroup of $G$ or subnormal in $G$. Arguing as above, $M_0$
cannot be subnormal in $G$. Then $M_0$ is a TI-subgroup of $G$ and so $M_0$ is also a TI-subgroup of $M$.
If $M_0$ is not normal in $M$, then $N_M(M_0)<M$. Let $M_1\leq M$ such that $N_M(M_0)$ is maximal in $M_1$, then
$N_M(M_0)\trianglelefteq M_1$. Take $x\in M_1\backslash N_M(M_0)$, one has $M_0\cap M_0^x=1$. Note that
$M_0^x\vartriangleleft(N_M(M_0))^x=N_M(M_0)$. Thus $M_0M_0^x=M_0\times M_0^x$. It follows that $(M_0\cap P)(M_0\cap P)^x=
(M_0\cap P)\times(M_0\cap P)^x$, this contradicts that $P$ is cyclic. Thus every non-trivial subgroup of $M$ of order
divisible by $p$ is normal in $M$. It follows that every non-trivial subgroup of $H$ is normal in $H$. Then $H$ is a
Dedekind group (see {\rm\cite[Theorem 5.3.7]{rob}}. Moreover, since every Sylow $t$-subgroup of $H$ is cyclic if $t>2$
and cyclic or a generalized quaternion group if $t=2$ by {\rm\cite[Theorem 10.5.6($ii)$]{rob}}), it follows that $H$ is
cyclic or a direct product of a quaternion group $Q_8$ and a cyclic group of odd order.

{\bf Subcase (2)}: Assume that $M$ has a maximal subgroup $M_2$ of $p'$-order and $M_2\ntrianglelefteq M$. Then
$M=P\rtimes M_2$. Note that $M=P\rtimes H$, we can assume $M_2=H$. By the maximality of $H$, one has that $P$ is
a minimal normal subgroup of $M$. Since $P$ is cyclic by {\rm\cite[Theorem 10.5.6($ii)$]{rob}}), it follows that
$P=Z_p$. Then $M=Z_p\rtimes H$. Arguing as in Subcase (1), one can get that $H$ is cyclic or a direct product
of a quaternion group $Q_8$ and a cyclic group of odd order.

In the following we prove the sufficiency part.

$(a)$ Suppose that $G$ is a group belonging to case (1), the proof is trivial.

$(b)$ Suppose that $G$ is a group belonging to case (2). Let $A$ be any subgroup of $G$ of even order. If $Z_q<A$,
then $A=A\cap (Z_q\rtimes\langle a\rangle)=Z_q\rtimes(A\cap\langle a\rangle)\trianglelefteq Z_q\rtimes\langle a\rangle=G$.
If $Z_q\nless A$, then $(q, |A|)=1$. Since $\langle a\rangle $ is a Hall-subgroup of $G$ and $G$ is obvious solvable,
we can assume $A\leq \langle a\rangle$ by {\rm\cite[Theorem 9.1.7]{rob}}. It is easy to see that $N_G(A)=\langle a\rangle$
since $\langle a\rangle$ is maximal in $G$. Then $A^g\cap A\leq \langle a\rangle^g\cap \langle a\rangle=1$ for each
$g\in G\backslash N_G(A)=G\backslash \langle a\rangle$, one has that $A$ is a TI-subgroup of $G$.

$(c)$ Suppose that $G$ is a group belonging to case (3). Let $A$ be any subgroup of $G$ of order divisible by $p$.
If ${Z_q}^r<A$, then $A=A\cap ({Z_q}^r\rtimes(P\times H))={Z_q}^r\rtimes(A\cap(P\times H))\trianglelefteq
{Z_q}^r\rtimes(P\times H)=G$. If ${Z_q}^r\nless A$, since $A\cap P\neq 1$ and every non-identity subgroup of
$P$ acts irreducibly on ${Z_q}^r$, one has that ${Z_q}^r\cap A=1$ and $P\times H$ is maximal in $G$. Then we
can assume $A\leq P\times H$. Note that $P\times H$ is a Dedekind-group and $A\ntrianglelefteq G$. One has
$N_G(A)=P\times H$. Therefore, $A^g\cap A\leq (P\times H)^g\cap (P\times H)=1$ for each $g\in G\backslash N_G(A)
=G\backslash(P\times H)$. It implies that $A$ is a TI-subgroup of $G$.

$(d)$ Suppose that $G$ is a group belonging to case (4). Let $A$ be any subgroup of $G$ of order divisible by $p$.
If ${Z_q}^r<A$, then $A=A\cap ({Z_q}^r\rtimes(Z_p\rtimes H))={Z_q}^r\rtimes(A\cap(Z_p\rtimes H))=
{Z_q}^r\rtimes(Z_p\rtimes(A\cap H))\trianglelefteq {Z_q}^r\rtimes(Z_p\rtimes H)$ since $H$ is a Dedekind-group.
If ${Z_q}^r\nless A$, arguing as above, one has that ${Z_q}^r\cap A=1$ and $Z_p\rtimes H$ is a maximal subgroup of
$G$. We can assume $A\leq Z_p\rtimes H$. Then $A=A\cap(Z_p\rtimes H)=Z_p\rtimes(A\cap H)\trianglelefteq Z_p\rtimes H$.
It follows that $N_G(A)=Z_p\rtimes H$. One has that $A^g\cap A\leq (Z_p\rtimes H)^g\cap (Z_p\rtimes H)=1$ for each
$g\in G\backslash N_G(A)=G\backslash (Z_p\rtimes H)$. It shows that $A$ is a TI-subgroup of $G$.

\hfill\qed

\section{Proof of Theorem~\ref{th001}}\label{th002}

\demo We only need to prove the necessity part.

Suppose that the theorem is false. Assume that $K$ is the largest subgroup of $G$ of order divisible by $p$ that is neither a
TI-subgroup of $G$ nor subnormal in $G$, then for any subgroup $L$ of $G$ satisfying $L>K$ we have that $L$ is a TI-subgroup of $G$
or subnormal in $G$. It follows that $K$ is not self-centralizing in $G$ by the hypothesis and then $K<N_G(K)$.

Considering the following subgroups series: $K<N_G(K)=K_1\leq N_G(K_1)=K_2\leq N_G(K_2)=K_3\leq N_G(K_3)=\cdots=K_{t-1}\leq
N_G(K_{t-1})=K_t\leq N_G(K_t)\leq\cdots$.

(1) Suppose that there exists a positive integer $t\geq 1$ such that $N_G(K_t)=G$. It follows that $K$ is subnormal in $G$, a
contradiction.

(2) Suppose that for any positive integer $t\geq 1$ we have $N_G(K_t)<G$. Note that $G$ is a finite group. It follows that there must
exist a positive integer $t\geq 1$ such that $K_t=N_G(K_t)<G$, which implies that $K_t$ is self-centralizing in $G$. It is obvious that
$K_t>K$. Arguing as in (1), $K_t$ cannot be subnormal in $G$. Then $K_t$ is a non-normal TI-subgroup of $G$ by the hypothesis. Moreover,
since $K_t=N_G(K_t)$, one has that $G$ is a Frobenius group with complement $K_t$. Let $N$ be the kernel of $G$, then $G=N\rtimes K_t$.

{\bf Claim I:} $K_t$ is either nilpotent or non-nilpotent and $K_t=Z_p\rtimes M$, where $Z_p\in\Syl_p(K_t)$, $M$ is nilpotent and $M$
has $p'$-order satisfying $[Z_p,M]\neq 1$.

Suppose that $K_t$ is non-nilpotent, then $K_t$ has at least one non-normal maximal subgroup. Let $M$ be a non-normal maximal subgroup
of $K_t$, then $N\rtimes M$ is a non-normal maximal subgroup of $G$. It follows that $N\rtimes M$ is self-centralizing in $G$ since
$N\rtimes M=N_G(N\rtimes M)$. Assume $p\mid|M|$, then $p\mid|N\rtimes M|$. By the hypothesis, $N\rtimes M$ is a TI-subgroup of $G$ or
subnormal in $G$. If $N\rtimes M\ntrianglelefteq G$, then $N\rtimes M$ is a non-normal TI-subgroup of $G$. However, since
$(N\rtimes M)^g\cap(N\rtimes M)=(N^g\rtimes M^g)\cap(N\rtimes M)=(N\rtimes M^g)\cap(N\rtimes M)\geq N>1$ for each $g\in G$, it follows
that $N\rtimes M\trianglelefteq G$, a contradiction. Therefore $N\rtimes M\trianglelefteq G$. One has $M=(N\cap K_t)M=(N\rtimes M)\cap
K_t\trianglelefteq K_t$, a contradiction. Thus $p\nmid|M|$. Let $P\in\Syl_p(K_t)$, arguing as above, one has $P\trianglelefteq K_t$.
Then $K_t=P\rtimes M$. Since $M$ is maximal in $K_t$, $P$ is a minimal normal subgroup of $K_t$, which implies that $P$ is an elementary
abelian group. Moreover, since $K_t$ is a complement of Frobenius group $G$, by {\rm\cite[Theorem 10.5.6(ii)]{rob}} one can get that
$P$ must be a cyclic group of prime order. Let $P=Z_p$, then $K_t=Z_p\rtimes M$. Moreover, arguing as above, one can obtain that every
maximal subgroup of $M$ is normal in $M$ and then $M$ is nilpotent.

{\bf Claim II:} $K\trianglelefteq K_t$.

Case $(i)$: Assume that $K_t$ is nilpotent. If $K_t$ is cyclic, then $K\trianglelefteq K_t$. Next assume that $K_t$ is non-cyclic.
Let $P\in\Syl_p(K_t)$.

If $p=2$, take an element $x$ of $P$ of order 2. Then $x$ induces a fixed-point-free automorphism of $N$. It follows that $N$ is an
abelian group of odd order and $g^x=g^{-1}$ for each non-identity element $g\in N$. Let $g$ be an element of $N$ of prime order,
then $\langle g\rangle\rtimes\langle x\rangle$ is also a Frobenius group of order divisible by $2$. Since $C_G(\langle g\rangle\rtimes
\langle x\rangle)\leq C_G(\langle g\rangle)\cap C_G(\langle x\rangle)=N\cap C_G(\langle x\rangle)=C_N(\langle x\rangle)=1<\langle g\rangle
\rtimes\langle x\rangle$, $\langle g\rangle\rtimes\langle x\rangle$ is self-centralizing in $G$. By the hypothesis, $\langle g\rangle\rtimes
\langle x\rangle$ is a TI-subgroup of $G$ or subnormal in $G$. It follows that $\langle g\rangle\rtimes\langle x\rangle$ is a TI-subgroup
of $N\rtimes\langle x\rangle$ or subnormal in $N\rtimes\langle x\rangle$. $(a)$ Assume that $\langle g\rangle\rtimes\langle x\rangle$
is a TI-subgroup of $N\rtimes\langle x\rangle$. Since $(\langle g\rangle\rtimes\langle x\rangle)^y\cap(\langle g\rangle\rtimes\langle x\rangle)=
(\langle g\rangle^y\rtimes\langle x\rangle^y)\cap(\langle g\rangle\rtimes\langle x\rangle)=(\langle g\rangle\rtimes\langle x\rangle^y)\cap
(\langle g\rangle\rtimes\langle x\rangle)\geq \langle g\rangle\neq 1$ for each $y\in N\rtimes\langle x\rangle$, one has $\langle g\rangle
\rtimes\langle x\rangle\trianglelefteq N\rtimes\langle x\rangle$. It follows that $N<\langle g\rangle\rtimes\langle x\rangle$ as
$\langle g\rangle\rtimes\langle x\rangle\nleq N$. Therefore $N=N\cap(\langle g\rangle\rtimes\langle x\rangle)=\langle g\rangle(N\cap
\langle x\rangle)=\langle g\rangle$. By N/C-theorem, $K_t\cong G/N=N_G(N)/C_G(N)\lesssim\Aut(N)$ is cyclic. It follows that $K_t$ is
cyclic and then $K\trianglelefteq K_t$. $(b)$ Assume that $\langle g\rangle\rtimes\langle x\rangle$ is subnormal in $N\rtimes
\langle x\rangle$. We claim $N\leq \langle g\rangle\rtimes\langle x\rangle$. Otherwise, assume $N\nleq \langle g\rangle\rtimes
\langle x\rangle$. Let $\langle g\rangle\rtimes\langle x\rangle=L_0\vartriangleleft L_1\vartriangleleft\cdots\vartriangleleft
L_{u-1}\vartriangleleft L_u=N\rtimes\langle x\rangle$ be a subnormal subgroups series, where $u\geq 1$. Since $N\nleq L_0$ and
$N\leq L_u$, there must exist a positive integer $v$ where $1\leq v\leq u$ such that $N\nleq L_{v-1}$ and $N\leq L_v$. Note that
$\langle g\rangle\rtimes\langle x\rangle\nleq N$, one has $N<L_v$. Then $L_v=N\rtimes(L_v\cap K_t)$, where $L_v\cap K_t>1$. It
follows that $L_v$ is a Frobenius group with complement $N$. Since $L_{v-1}\vartriangleleft L_v$ and $N\nleq L_{v-1}$, one has
$L_{v-1}<N$, which implies that $\langle g\rangle\rtimes\langle x\rangle=L_0<N$, a contradiction. Thus $N\leq \langle g\rangle
\rtimes\langle x\rangle$. Arguing as above, one has $N=\langle g\rangle$ and then $K_t$ is cyclic. One also has $K\trianglelefteq K_t$.

If $p>2$, let $P_2\in\Syl_2(K_t)$. Since $K_t$ is non-cyclic, $P_2$ is a generalized quaternion group of order $4n$ by
{\rm\cite[Theorem 10.5.6(ii)]{rob}}, where $n\geq 2$. Then $K_t=R\times P_2$, where $R$ is a cyclic group of order divisible by $p$. Let
$P_2=\langle a,b\mid a^n=b^2, a^{2n}=1, b^{-1}ab=a^{-1}\rangle$, one has $Z(P_2)=a^n=b^2$. For any non-identity subgroup $Q$
of $P_2$, it is easy to see that $C_G(R\times Q)=R\times Z(P_2)\leq R\times Q$ and then $R\times Q$ is self-centralizing in $G$.
By the hypothesis, $R\times Q$ is a TI-subgroup of $G$ or subnormal in $G$. Arguing as above, $R\times Q$ cannnot be subnormal in
$G$. Then $R\times Q$ is a TI-subgroup of $G$. It follows that $R\times Q$ is a TI-subgroup of $K_t$. If $R\times Q\ntrianglelefteq K_t$,
then there exists $g\in K_t$ such that $(R\times Q)^g\cap (R\times Q)=1$. However, $(R\times Q)^g\cap (R\times Q)=(R^g\times Q^g)\cap
(R\times Q)=(R\times Q^g)\cap(R\times Q)\geq R\neq 1$, a contradiction. Therefore $R\times Q\trianglelefteq K_t$. It follows that
$Q\trianglelefteq P_2$. By the choice of $Q$, one has that $P_2$ is a Dedekind-group. Then $K_t=R\times P_2$ is a Dedekind-group.
It follows that $K\trianglelefteq K_t$.

Case $(ii)$: Assume $K_t=Z_p\rtimes M$, where $M$ is nilpotent, $p\nmid|M|$ and $[Z_p,M]\neq 1$. Since $p\mid|K|$, one has
$Z_p\leq K$. Then $K=K\cap(Z_p\rtimes M)=Z_p\rtimes(K\cap M)$. If $M$ is cyclic, it is obvious that $Z_p\rtimes(K\cap M)\trianglelefteq Z_p\rtimes M$,
that is $K\trianglelefteq K_t$. If $M$ is non-cyclic, then $M=T\times P_2$, where $T$ is cyclic and $P_2$ is a generalized quaternion
group. For any non-identity subgroup $Q$ of $P_2$, $Z_p\rtimes(T\times Q)$ is a subgroup of $K_t$ of order divisible by $p$.
Moreover, it is clear that $C_G(Z_p\rtimes(T\times Q))\leq Z_p\rtimes(T\times Q)$, that is $Z_p\rtimes(T\times Q)$ is self-centralizing
in $G$. By the hypothesis, $Z_p\rtimes(T\times Q)$ is a TI-subgroup of $G$ or subnormal in $G$. Arguing as in case $(i)$, one can get
that $P_2$ must be a quaternion group of order 8 and then $M=T\times P_2$ is a Dedekind-group. It follows that
$K=Z_p\rtimes(K\cap M)\trianglelefteq Z_p\rtimes M=K_t$.

By claim II, one has $K_t\leq N_G(K)$. It follows that $N_G(K)=K_t$ as $N_G(K)=K_1\leq K_t$. Since $K^g\cap K\leq {K_t}^g\cap K_t=1$
for each $g\in G\backslash N_G(K)=G\backslash K_t=G\backslash N_G(K_t)$, one has that $K$ is a TI-subgroup of $G$, a
contradiction.

It implies that our assumption is not true and so every subgroup of $G$ is a TI-subgroup or subnormal or has $p'$-order.
\hfill\qed

\end{document}